\documentclass[11pt]{amsart}
\usepackage{graphicx}

\usepackage{amsmath}
\usepackage{amsthm}
\usepackage{amssymb}
\usepackage{enumerate}
\usepackage{a4wide}

\def\AC{\mathcal{AC}}
\def\pf{\begin{proof}}
\def\pfk{\end{proof}}

\newtheorem{theorem}{Theorem}[section]

\newtheorem{lemma}[theorem]{Lemma}
\newtheorem{corollary}[theorem]{Corollary}
\newtheorem{remark}[theorem]{Remark}
\newtheorem{proposition}[theorem]{Proposition}
\newtheorem{observation}[theorem]{Observation}
\newtheorem{example}[theorem]{Example}

\begin{document}
\title{Abelian Complexity of Infinite Words Associated with Quadratic Parry Numbers}

\author[Balkov\'a]{L{\!'}ubom\'{\i}ra~Balkov\'a}
\address[Balkov\'a]{Department of Mathematics, FNSPE Czech Technical University in
Prague, Trojanova 13, 120 00 Praha 2, Czech Republic}
\email{lubomira.balkova@fjfi.cvut.cz}

\author[B\v rinda]{Karel B\v rinda}
\address[B\v rinda]{Department of Mathematics, FNSPE Czech Technical University in
Prague, Trojanova 13, 120 00 Praha 2, Czech Republic}
\email{karel.brinda@fjfi.cvut.cz}

\author[Turek]{Ond\v rej Turek}
\address[Turek]{Laboratory of Physics, Kochi University of Technology,
Tosa Yamada, Kochi 782-8502, Japan}
\email{ondrej.turek@kochi-tech.ac.jp}
\keywords{Abelian complexity, quadratic Parry numbers, binary infinite words}
\subjclass[2010]{68R15}
\date{\today}

\begin{abstract}
We derive an explicit formula for the Abelian complexity of infinite words associated with quadratic Parry numbers.

\end{abstract}
\maketitle
\section{Introduction}
Abelian complexity is a~property of infinite words that has been examined for the first time by Coven and Hedlund in~\cite{CoHe}, where they have revealed that it could serve as an alternative way of characterization of periodic words and of Sturmian words. They have shown that an infinite word is periodic if and only if its Abelian complexity satisfies $\AC(n) = 1$ for large enough $n$, and proven that an aperiodic binary infinite word is Sturmian if and only if its Abelian complexity is identically equal to $2$.

However, the notion ``Abelian complexity'' itself comes from the paper~\cite{RSZ} that in a~sense initiated a~general study of the Abelian complexity of infinite words over finite alphabets. It is noteworthy that besides Sturmian words, the Abelian complexity is known still for only few infinite words.
To answer affirmatively the question of G. Rauzy: ``Does there exist an infinite word $\mathbf u$ whose Abelian complexity $\AC(n)=3$ for all $n \in \mathbb N$?'',
the authors of~\cite{RSZ} have found two classes of words with $\AC(n)=3$ for all $n \in \mathbb N$.
Ibidem, the Abelian complexity of the Thue-Morse word has been described and the class of words having the same Abelian complexity as the Thue-Morse word has been characterized.
The Abelian complexity of the Tribonacci word has been studied in~\cite{RSZ2} and it has been shown that $\AC(n) \in \{3,4,5,6,7\}$ for all $n \in \mathbb N$, each of these five values is assumed,
and the values $3$ and $7$ are attained infinitely many times.
Moreover, the set $\{n \in \mathbb N \bigm | \AC(n)=3\}$ has been described. One of the coauthors of this paper has provided in~\cite{Tu2} an optimal upper bound on the Abelian complexity
for ternary infinite words associated with cubic Pisot numbers, roots of polynomials $x^3-px^2-x+1$ with $p>1$, and has shown that the upper bound is reached infinitely many times.

In this paper, we will determine the Abelian complexity of infinite words associated with quadratic Parry numbers.
\section{Preliminaries}\label{Preliminaries}
An {\em alphabet} $\mathcal A$ is a~finite set of symbols called {\em
letters}; in the paper we will use only the binary alphabet ${\mathcal A}=\{A,B\}$.
Any finite sequence of letters from ${\mathcal A}$ is called a~{\em word} (over ${\mathcal A}$). The set
$\mathcal A^{*}$ of all finite words (including the empty word
$\varepsilon$) provided with the operation of concatenation is
a~free monoid. If $w = w_0w_1w_2\cdots w_{n-1}$ is a finite word over ${\mathcal A}$, we denote its length by $|w| = n$, and use the symbols $|w|_A, |w|_B$ for the number of occurrences of the letters $A$ and $B$ in $w$, respectively. The {\em mirror image} of the word $w = w_0w_1 \cdots w_{n-1}$ is defined as $\overline{w} = w_{n-1}\cdots w_1w_0$.
We will deal also with infinite sequences of letters from ${\mathcal A}$, called {\em infinite words} ${\bf u}=u_0u_1u_2\cdots$ (over ${\mathcal A}$).
A~finite word $w$ is called a~{\em factor} of
the word ${\bf u}$ (${\bf u}$ being finite or infinite) if there exist a~finite word
$p$ and a~word $s$ (finite or infinite) such that
${\bf u}=pws$. We say that the word
$w$ is a~prefix of ${\bf u}$ if $p = \varepsilon$, and a~suffix of ${\bf u}$ if $s = \varepsilon$. A prefix $p$ of ${\bf u}$ is said to be proper if $p\neq{\bf u}$, a proper suffix is defined in a similar way. A~concatenation of $k$ words $w$
is denoted by $w^k$, a~concatenation of infinitely many finite words $w$ by $w^\omega$. An infinite word ${\bf u}$ is said to be {\em
eventually periodic} if there exist words $v,w$ such that ${\bf u}=v
w^{\omega}$. A~word which is not eventually periodic is called
{\em aperiodic}.

Let us associate with every factor $w$ of a~binary infinite word ${\bf u}$ its {\em Parikh vector} $\Psi(w)=(|w|_A, |w|_B)$.
The {\em Abelian complexity} is the map ${\mathcal AC}: \mathbb N \to \mathbb N$ defined by
$${\mathcal AC}(n)=\#\{\Psi(w) \, | \, w \ \text{is a factor of ${\bf u}$ of length $n$}\}.$$
(Throughout this paper $\mathbb N$ stands for the set of positive integers and ${\mathbb N}_0 = \mathbb N \cup \{0\}$.)
When considering two consecutive factors of the same length of an infinite word, then
the corresponding entries of their Parikh vectors can differ at most by one.
As a~consequence, we obtain the following observation.
\begin{observation}\label{Parikh_by_one}
Let $(k,n-k)$ and $(k',n-k')$, where $k < k'$, be Parikh vectors of factors of a~binary infinite word ${\bf u}$.
Then $(k'',n-k'')$ is the Parikh vector of a factor of $\bf u$ for any $k'' \in \mathbb N, \ k \leq k'' \leq k'$.
\end{observation}

Using the previous observation, we can express $\AC(n)$ as follows.
\begin{proposition}\label{AC(n)}
Let ${\bf u}$ be a~binary infinite word. Then for all $n \in \mathbb N$
$${\AC}(n)=1+\max\{\,|\,|v|_A-|w|_A| \ \bigm | \ v, \ w \ \text{factors of $\bf u$ of length $n$}\}.$$
\end{proposition}

Proposition~\ref{AC(n)} reflects a~close relation between the Abelian complexity and the balance property of binary infinite words.
A~binary infinite word ${\bf u}$ is said to be $C$-{\em balanced} if
for any pair of factors $w, v$ of ${\bf u}$, with $|w| = |v|$, it holds $||w|_A - |v|_A|\leq C$. Let us remark that Sturmian words can be
defined as aperiodic $1$-balanced binary infinite words~\cite{MoHe}.
Combining Proposition~\ref{AC(n)} with the definition of binary $C$-balanced words, we can draw the following conclusion.
\begin{corollary}\label{balance_AC}
Let ${\bf u}$ be a binary infinite word. Then its optimal balance bound is $C$ if and only if the maximum of its Abelian complexity is $C+1$.
\end{corollary}
\section{Infinite words associated with quadratic Parry numbers}\label{Ubeta}
The main aim of this paper is to determine the Abelian complexity of infinite words associated with quadratic Parry numbers.
These words are usually defined via the so-called $\beta$-integers corresponding to Parry numbers $\beta$, for details see~\cite{BPT}. However, here we skip the procedure for the sake of brevity and introduce them directly using the notion of morphism. A~mapping $\varphi: \{A,B\}^{*} \to \{A,B\}^{*}$ is
called a~{\em morphism} if $\varphi(vw)=\varphi(v)\varphi(w)$ for
all $v,w \in \{A,B\}^{*}$. Obviously, a morphism is uniquely
determined by $\varphi(A)$ and $\varphi(B)$. The action of
the morphism $\varphi$ can be naturally extended to
infinite words by
$$
\varphi(u_0u_1u_2\cdots):=\varphi(u_0)\varphi(u_1)\varphi(u_2)\cdots
$$
An infinite word ${\bf u}$
such that $\varphi({\bf u})={\bf u}$ is called a~{\em fixed point} of the
morphism $\varphi$.
To any morphism $\varphi$ on $\{A,B\}$, one can associate the {\it
incidence matrix} $M_\varphi$ defined by
$$
M_\varphi=\left(\begin{matrix}|\varphi(A)|_A & |\varphi(A)|_B\\
|\varphi(B)|_A & |\varphi(B)|_B \end{matrix}\right).
$$
It follows immediately from the definition of $M_\varphi$ that for any word $u \in \{A,B\}^*$
\begin{equation}\label{matice}
(|\varphi(u)|_A,|\varphi(u)|_B)=(|u|_A, |u|_B)M_\varphi\,.
\end{equation}

It has been shown in~\cite{Fa} that infinite words ${\bf u}_\beta$ associated with quadratic Parry numbers $\beta$ are in fact fixed points of certain morphisms. These morphisms are of two types:

$\bullet$\quad The class of morphisms providing infinite words associated with {\em quadratic simple Parry numbers} is given by
\begin{equation}\label{simple}
\varphi: A\mapsto A^p B, \quad B\mapsto A^q \quad \text{for $p,q\in\mathbb N$, $p\geq q$.}
\end{equation}
The corresponding incidence matrix is
$$
M_\varphi=\left(\begin{array}{cc}p&1\\q&0\end{array}\right)\,.
$$
The infinite word ${\bf u}_\beta$ is the only fixed point of this morphism, i.e., ${\bf u}_\beta=\varphi({\bf u}_\beta)$,
and it can be obtained when iterating $\varphi$ on the letter $A$ infinitely many times. Formally, ${\bf u}_\beta=\lim_{n \to \infty}\varphi^n(A)$ taken with respect to the product topology.

$\bullet$\quad The class of morphisms providing infinite words associated with {\em quadratic non-simple Parry numbers} is given by
\begin{equation}\label{nonsimple}
\varphi: A\mapsto A^p B, \quad B\mapsto A^q B \quad \text{for $p,q\in\mathbb N$, $p>q$.}
\end{equation}
The associated incidence matrix is
$$
M_\varphi=\left(\begin{array}{cc}p&1\\q&1\end{array}\right)\,.
$$
Again, ${\bf u}_\beta=\lim_{n \to \infty}\varphi^n(A)$ is the only fixed point of $\varphi$.

Throughout the whole paper, the symbol $\varphi$ denotes a morphism associated with a quadratic (simple or non-simple) Parry number
and ${\bf u}_\beta$ is exclusively used for the word $\lim_{n\to\infty}\varphi^n(A)$.

For infinite words associated with quadratic Parry numbers, the optimal balance bounds are already known~\cite{Tu1, BPT}.
By Corollary~\ref{balance_AC}, the maximum of the Abelian complexity is immediately determined.
\begin{corollary}\label{maxAC_Parry}
The Abelian complexity of ${\bf u}_\beta$ satisfies
\begin{itemize}
\item $\max\{\AC(n) \bigm | n \in \mathbb N\}=2+\lfloor \frac{p-1}{p+1-q}\rfloor$ in the case of simple Parry numbers,
\item $\max\{\AC(n) \bigm | n \in \mathbb N\}=1+\lceil \frac{p-1}{q}\rceil$ in the case of non-simple Parry numbers.
\end{itemize}
\end{corollary}

The target of this paper is to derive an explicit formula for $\AC(n)$ for all $n\in\mathbb N$.

\begin{remark}\label{AC_Sturmian}
It is seen from Corollary~\ref{maxAC_Parry} that for a~simple Parry number $\beta$ with $q=1$ or for a~non-simple Parry number $\beta$ with $p=q+1$, the Abelian complexity satisfies $\AC(n)=2$ for all $n \in \mathbb N$.
This result follows also from~\cite{FrGaKr}, where it has been shown that for such $\beta$ the infinite word ${\mathbf u}_\beta$ is Sturmian, thus its Abelian complexity is constant and equal to $2$ by~\cite{CoHe}. Moreover, these cases are the only ones among all Parry numbers $\beta$ (not only the quadratic ones) for which the infinite word ${\mathbf u}_\beta$ is Sturmian, see~\cite{FrMaPe}.
\end{remark}

For the calculation of the Abelian complexity, the following simple observation deduced from~\eqref{matice} will be important.
\begin{observation}\label{phi^j(A)}
For all $n\in\mathbb N$, it holds
$$\Psi\left(\varphi^n(A)\right)=(1,0)M_\varphi^n\,, \quad \quad |\varphi^n(A)|_A=(1,0)M_\varphi^n\left(\begin{matrix}1\\0\end{matrix}\right), \quad \quad |\varphi^n(A)|_B=(1,0)M_\varphi^n\left(\begin{matrix}0\\1\end{matrix}\right).$$
\end{observation}
Let us define a~strictly increasing sequence $U=(U_n)_{n=0}^{\infty}$ by
$U_n=|\varphi^n(A)|$, i.e.,
\begin{equation}\label{Un}
U_n=(1, 0)M_\varphi^n \left(\begin{array}{c}1\\1\end{array}\right)\,.
\end{equation}
Furthermore, let us denote by
$\langle n\rangle_U$ the normal $U$-representation of $n$; for details on $U$-representations see Lothaire~\cite{Lo}.
Recall that the normal $U$-representation of $n \in \mathbb N$ is equal to
\begin{equation}\label{Urepr}
\langle n \rangle_U=(d_N, d_{N-1}, \ldots, d_1, d_0)
\end{equation}
if $n=\sum_{j=0}^N d_j U_j$ and the representation is obtained by the following algorithm:
\begin{enumerate}
\item Find $N\in\mathbb N_0$ such that $n<U_{N+1}$.
\item Put $x_{N}:=n$.
\item For $j=N, N-1,\ldots,1,0$ do $$d_j:=\left\lfloor\frac{x_j}{U_j}\right\rfloor \quad \quad \text{and}\quad \quad x_{j-1}:=x_j-d_jU_j.$$
\end{enumerate}
\begin{remark}
Note that $N$ can be chosen as any integer such that $n<U_{N+1}$, i.e., not necessarily the smallest one satisfying $U_N\leq n<U_{N+1}$. As a result, the normal $U$-representation of $n$ can start with a block of zeros, and the representation $(0,0,...,0,d_N,d_{N-1},...,d_0)$ is equivalent to $(d_N,d_{N-1},...,d_0)$.
\end{remark}
Since the inequality
$U_{N+1}=|(\varphi^N(A))^p\varphi^N(B)|<|(\varphi^N(A))^{p+1}|=(p+1)U_N$
holds for all $N \in \mathbb N_0$, one can see that the coefficients in normal $U$-representations are less than or equal to~$p$.

The following proposition is a~particular case of a~general theorem concerning morphisms associated with Parry numbers provided in~\cite{Fa}.
\begin{proposition}\label{struktura pref. u}
Let $u$ be the prefix of ${\mathbf u}_\beta$ of length $n$, $n \in \mathbb N$. If $\langle n\rangle_{U}=(d_{N},d_{N-1},\ldots,d_1,d_0)$, then
$u = \left(\varphi^{N}(A)\right)^{d_{N}}\left(\varphi^{N-1}(A)\right)^{d_{N-1}}\cdots\left(\varphi(A)\right)^{d_1}A^{d_0}$, where the zero-power of a~word $w$ is defined as the empty word, i.e., $w^0=\varepsilon$.
\end{proposition}
\section{The main idea}
Let us open the calculation of $\AC(n)$ by this short section in which we will describe the key idea of our technique.

It turns out that for any simple or non-simple quadratic Parry number $\beta$, one can find two infinite words, let us denote them by ${\bf v}$ and ${\bf w}$, that have the following properties.
\begin{itemize}
\item For any finite prefix $\hat{v}$ of ${\bf v}$ it holds:
\begin{enumerate}
\item $\hat{v}$ is a factor of ${\bf u}_\beta$.
\item For any factor $u$ of ${\bf u}_\beta$,
$$
|u|=|\hat{v}| \quad\Rightarrow\quad |u|_B\geq|\hat{v}|_B\,.
$$
\end{enumerate}
In other words, the number of occurrences of letter $B$ in factors of ${\bf u}_\beta$ of a given length $n$ attains its \emph{minimum} in the prefix of ${\bf v}$ of length $n$.
\item For any finite prefix $\hat{w}$ of ${\bf w}$ it holds:
\begin{enumerate}
\item $\hat{w}$ is a factor of ${\bf u}_\beta$.
\item For any factor $u$ of ${\bf u}_\beta$,
$$
|u|=|\hat{w}| \quad\Rightarrow\quad |u|_B\leq|\hat{w}|_B\,.
$$
\end{enumerate}
In other words, the number of occurrences of letter $B$ in factors of ${\bf u}_\beta$ of a given length $n$ attains its \emph{maximum} in the prefix of ${\bf w}$ of length $n$.
\end{itemize}

Note that ${\bf v}$ and ${\bf w}$ are unique. These words ${\bf v}$ and ${\bf w}$ play an essential role in computation of the Abelian complexity.
\begin{proposition}\label{AC=rozdil}
Let ${\bf v}$ and ${\bf w}$ be the infinite words defined above.
Then the Abelian complexity of ${\bf u}_\beta$ can be expressed for all $n \in \mathbb N$ by the formula
\begin{equation}\label{formule}
\AC(n)=1+|\hat{w}|_B-|\hat{v}|_B\,,
\end{equation}
where $\hat{w}$ and $\hat{v}$ are prefixes of ${\bf w}$ and ${\bf v}$, respectively, of the same length $n$.
\end{proposition}

\pf
Using properties of $\bf v$ and $\bf w$, it follows that $$\max\{\left||v|_A-|w|_A\right| \bigm | \text{ $u, w$ factors of ${\bf u}_\beta$ of length $n$}\} = |\hat v|_A-|\hat w|_A=|\hat w|_B-|\hat v|_B,$$
where $\hat v, \ \hat w$ are prefixes of length $n$ of $\bf v$ and $\bf w$, respectively.
The statement is then a~direct consequence of Proposition~\ref{AC(n)}.
\pfk

For an application of the formula~\eqref{formule}, it is necessary to find the words ${\bf v}$ and ${\bf w}$ and their structure which will be done in the rest of the paper, separately for the simple and the non-simple Parry case.
\section{$\AC(n)$ in the non-simple Parry case}\label{Kap.NonSimple}
At first we will determine the Abelian complexity of infinite words ${\bf u}_\beta$ associated with quadratic non-simple Parry numbers, which are fixed points of the morphisms of the type
$$
\varphi(A)=A^pB\,,\quad \varphi(B)=A^qB\,,\qquad p>q\geq1\,.
$$
Let us recall that according to Corollary~\ref{maxAC_Parry}, the maximum of the Abelian complexity of these words equals $1+\lceil\frac{p-1}{q}\rceil$. Moreover, by Remark~\ref{AC_Sturmian}, in the special case $p=q+1$ it holds $\AC(n)=2$ for all $n\in\mathbb N$.
The \emph{non-simple} case -- despite its name -- is slightly easier to treat than the \emph{simple} one,
because the words ${\bf v}$ and ${\bf w}$ have been already found in the paper~\cite{BPT}.
Let us just adopt them from there:
$$
{\bf v}={\bf u}_\beta\,,
$$
$$
{\bf w} =  \lim_{n\to\infty}w^{(n)}\,,
$$
where $w^{(0)} = B \ \text{and} \ w^{(n)} = B\varphi(w^{(n-1)}) \ \text{for $n \in \mathbb N$}\,.$

In order to determine $|\hat v|_B$ and $|\hat w|_B$ for any prefix $\hat v$ of ${\mathbf v}={\mathbf u}_\beta$ and $\hat w$ of ${\mathbf w}$,
we will use the following observation, which is a~direct consequence of the morphism form.
\begin{observation}\label{Phi_B}
For any $u\in{\mathcal A}^*$ it holds $|\varphi(u)|_B=|u|$. In particular, $|\varphi^j(A)|_B=U_{j-1}$ for all $j\in\mathbb N$.
\end{observation}

With the previous observation and Proposition~\ref{struktura pref. u} in hand, we can compute $|\hat v|_B$ for any prefix $\hat v$ of ${\mathbf v}={\mathbf u}_\beta$.
\begin{lemma}\label{pocetBve_v}
Let $n \in \mathbb N$ with $\langle n\rangle_{U}=(d_{N},d_{N-1},\ldots,d_1,d_0)$ and let $\hat v$ be the prefix of ${\mathbf v}=\textbf{u}_\beta$ of length $n$. Then
$$
|\hat v|_B=\sum_{j=1}^{N}d_j U_{j-1}\,.
$$
\end{lemma}

The following lemma will play an essential role in calculation of $|\hat w|_B$ for any prefix $\hat w$ of~${\mathbf w}$.

\begin{lemma}\label{lemmatkoow}
For any $N \in {\mathbb N}_0$, the word $w^{(N)}$ has the following two properties:
\begin{enumerate}
  \item $w^{(N)}$ is a palindrome,
  \item $w^{(N)}$ is a proper suffix of $\varphi^{N+1}(A)$.
\end{enumerate}
\end{lemma}
\pf
Let us proceed by induction on $N$.

\begin{enumerate}
\item $w^{(0)} = B$ is a palindrome. Let $N \geq 0$ and let us suppose that $w^{(N)}$ is a palindrome. We will prove
that $w^{(N+1)}$ is a palindrome, too.
We have $w^{(N+1)} = B\varphi(w^{(N)}) = B\varphi(\overline{w^{(N)}})$.
Denote $w^{(N)}=w_0w_1\cdots w_n$, where $w_i$ are letters. Since for every $a \in {\mathcal A}$ it holds $B\varphi(a)=\overline{\varphi(a)}B$,
we deduce $B\varphi(\overline{w^{(N)}})=B\varphi(w_n)\cdots \varphi(w_1)\varphi(w_0)=$\\
$\overline{\varphi(w_n)}\cdots \overline{\varphi(w_1)}\ \overline{\varphi(w_0)}B=\overline{\varphi(w^{(N)})}B=\overline{B\varphi(w^{(N)})}=\overline{w^{(N+1)}}.$

\item $w^{(0)} = B$ is a~proper suffix of $\varphi(A)$.
Let $N \geq 1$ and let us assume $\varphi^N(A) = uw^{(N-1)}$, where $u \not =\varepsilon$.
Thus $\varphi^{N+1}(A) = \varphi(u)\varphi(w^{(N-1)}).$
Since the last letter of $\varphi(u)$ is $B$, we have proven that $w^{(N)}=B\varphi(w^{(N-1)})$ is a~proper suffix of $\varphi^{N+1}(A)$.
\end{enumerate}
\pfk

We have prepared everything for determining $|\hat w|_B$ for any prefix $\hat w$ of ${\mathbf w}$.
\begin{lemma}\label{pocetBve_w}
Let $n \in \mathbb N$ and let $\hat w$ be the prefix of ${\mathbf w}$ of length $n$. Take $k$ arbitrary such that $n \leq |w^{(k)}|$
and denote $\langle U_{k+1}-n\rangle_U=(e_{k}, e_{k-1},\ldots , e_1, e_0)$.
Then $$|\hat w|_B=U_k-\sum_{j=1}^{k}e_jU_{j-1}\,.$$
\end{lemma}
\pf
According to the choice of $k$, the word
$\hat w$ is a~prefix of $w^{(k)}$, and using the first statement of Lemma~\ref{lemmatkoow}, $\overline{\hat w}$ is a~suffix $w^{(k)}$.
By the second statement of Lemma~\ref{lemmatkoow}, we have $\varphi^{k+1}(A)=u\overline{\hat w}$, where $u$ is a prefix of ${\bf u}_\beta$, $u \not =\varepsilon$.
Therefore $|\hat w|_B=|\varphi^{k+1}(A)|_B-|u|_B$.
We obtain $|\varphi^{k+1}(A)|_B=U_{k}$ from Observation~\ref{Phi_B}, and
since $\langle |u|\rangle_U=\langle U_{k+1}-n\rangle_U=(e_{k},e_{k-1}, \ldots, e_1,e_0)$ by assumption,
it holds $|u|_B=\sum_{j=1}^{k}e_jU_{j-1}$ according to Lemma~\ref{pocetBve_v}.
\pfk

As an immediate consequence of Proposition~\ref{struktura pref. u} and Lemmas~\ref{pocetBve_v} and~\ref{pocetBve_w}, we get the main theorem.
\begin{theorem}\label{ACnon_simple}
Let ${\bf u}_\beta$ be the fixed point of the morphism $\varphi$ defined in~\eqref{nonsimple}.
Let $(U_n)_{n=0}^{\infty}$ be the sequence defined in~\eqref{Un}.
The Abelian complexity of the infinite word ${\bf u}_\beta$ is given for all $n \in \mathbb N$ by the formula
$$
\AC(n)= 1 + U_k - \sum_{j = 1}^{k}(d_j + e_j)U_{j-1}  \,,
$$
where
\begin{itemize}
\item $k$ is arbitrary such that $n \leq |w^{(k)}|$,
\item $\langle n\rangle_{U}=(d_{k},d_{k-1},\ldots,d_1,d_0)$,
\item $\langle U_{k+1} - n\rangle_{U}=(e_{k},e_{k-1},\ldots,e_1,e_0)$.
\end{itemize}
\end{theorem}

In the end of this chapter, let us estimate the minimal index $k$ from Theorem~\ref{ACnon_simple}.
\begin{lemma}\label{sendvic}
Let $n\in\mathbb N$ and let the number $N\in\mathbb N_0$ satisfy $U_{N}\leq n<U_{N+1}$. Then
$$
\left|w^{(N)}\right|\leq n<\left|w^{(N+2)}\right|\,.
$$
\end{lemma}
\pf
Let us recall that for all $N\in\mathbb N_0$, $w^{(N)}=B\varphi(B)\varphi^2(B)\cdots\varphi^{N}(B)$.
In order to prove the first inequality ($\leq$), it suffices to show that $|w^{(N)}|\leq U_N$. This can be done by induction on $N$, using the identities $\varphi^{N+1}(A)=\left(\varphi^{N}(A)\right)^p\varphi^{N}(B)$, $\varphi^{N+1}(B)=\left(\varphi^{N}(A)\right)^q\varphi^{N}(B)$ and the inequality $p > q$.
As for the second inequality ($<$), $n<U_{N+1}$ implies
$$
n<|\varphi^{N+1}(A)|<|\left(\varphi^{N+1}(A)\right)^q\varphi^{N+1}(B)|=|\varphi^{N+2}(B)|<\left|B\varphi(B)\cdots\varphi^{N+2}(B)\right|=|w^{(N+2)}|\,.
$$
\pfk
Lemma~\ref{sendvic} tells us that if $n\in\mathbb N$ satisfies $U_N\leq n<U_{N+1}$, then the minimal index $k$ from Theorem~\ref{ACnon_simple} is either $N+1$ or $N+2$. Hence the choice $k=N+2$ works universally.

\begin{example}
Let $p = 3$ and $q = 1$. Let us calculate $\AC(n)$ for $n=7$.
Using \eqref{Un}, we obtain $(U_n)_{n=0}^{\infty} = (1, 4, 14, 48, 164, \ldots)$.
Applying Lemma~\ref{sendvic}, we can put $k=3$ in Theorem~\ref{ACnon_simple}.
Now we have to get the normal $U$-representations:
$$\langle7\rangle_U = (0,0,1,3) = (d_3,d_2,d_1,d_0)
\hspace{0.2cm} \mathrm{and} \hspace{0.2cm}
\langle U_4 - 7 \rangle_U = \langle 157 \rangle_U = (3,0,3,1) = (e_3,e_2,e_1,e_0).$$
Finally, using Theorem~\ref{ACnon_simple}, we find $\AC (7) = 1 + U_3 - \sum_{j=1}^{3}(d_j + e_j)U_{j-1} = 3$.

Since $n$ is small, we can illustrate the situation by observing the prefix $\varphi^3(A)$ of ${\mathbf u}_\beta$:
$$ \varphi^{3}(A) = \underbrace{AAABAAA}_{\hat{v}} BAAABAB AAABAAABAAABAB
 AAABAAABAAA   \overbrace{BA\underbrace{B AAA   \overbrace{  B A\overbrace{B}^{  w^{(0)}  }  }^{w^{(1)}}     }_{\overline{\hat{w}}}}^{w^{(2)}} $$
We see that $ |\hat{w}|_B - |\hat{v}|_B = 2  $ and therefore indeed $\AC(7) = 3$. (Moreover, one can notice that in this case, it is sufficient to put $k=2$ in Theorem~\ref{ACnon_simple}.)
\end{example}
\section{$\AC(n)$ in the simple Parry case}\label{ACnSimple}
In this section we study the Abelian complexity of infinite words associated with quadratic simple Parry numbers, i.e., we deal with the fixed point ${\mathbf u}_\beta$ of the morphism $\varphi$ given by
$$
\varphi(A)=A^pB\,,\quad \varphi(B)=A^q\,,\qquad p\geq q\geq1\,.
$$

\subsection{Case $q=1$}

When $\beta$ is a quadratic simple Parry number, the calculation of $\AC(n)$ turns
out to be more effective when the cases $q=1$ and $q>1$ are treated separately.
Let us start with $q=1$. As we have explained in Remark~\ref{AC_Sturmian}, in this case it holds
$\AC(n)=2$ for all $n\in\mathbb N$.

\subsection{Case $q>1$}
From now on, let us assume that $q>1$. By Corollary~\ref{maxAC_Parry}, the maximum of the Abelian complexity is $2+\lfloor\frac{p-1}{p+1-q}\rfloor$.
In order to find an explicit expression for $\AC(n)$ for all $n\in\mathbb N$, we have to describe the infinite words $\mathbf v$ and $\mathbf w$ from Proposition~\ref{AC=rozdil}. These words, as it will be proven in Proposition~\ref{w v} below, are:
$$
{\bf w} =  \lim_{n\to\infty}w^{(n)}\,,
$$
where

$$
\begin{array}{rl}
w^{(0)}=&B \\
A^pw^{(n)}=&\varphi^2(w^{(n-1)})  \quad \text{for all $n \in \mathbb N$}\,,
\end{array}
$$
and
$$
{\bf v}= \lim_{n\to\infty}v^{(n)}\,,
$$
where $v^{(n)}=\varphi(w^{(n)}) \quad \text{for all $n \in \mathbb N_0$}$.

To demonstrate how $\mathbf v$ and $\mathbf w$ are constructed, here are their prefixes written explicitly:
\begin{equation}\label{defwvs}
\begin{array}{rcll}
w^{(n)}&=&B\varphi(A^{q-1})\varphi^3(A^{q-1})\varphi^5(A^{q-1})\cdots\varphi^{2n-3}(A^{q-1})\varphi^{2n-1}(A^{q-1})\,, \\
v^{(n)}&=&A^q\varphi^2(A^{q-1})\varphi^4(A^{q-1})\varphi^6(A^{q-1})\cdots\varphi^{2(n-1)}(A^{q-1})\varphi^{2n}(A^{q-1})\,.
\end{array}
\end{equation}
Both $w^{(n)}$ and $v^{(n)}$ are factors of ${\mathbf u}_\beta$: it can be easily proven by induction on $n$ that $w^{(n)}$ is a~suffix of $\varphi^{2n}(B)$, and $v^{(n)}$ is the image by $\varphi$ of $w^{(n)}$.

We observe from the construction of the words $\mathbf v$ and $\mathbf w$ that they are related via $\varphi$:
\begin{observation}\label{w v obs}
The inifinite words $\mathbf v$ and $\mathbf w$ satisfy the relations
$$
{\mathbf v}=\varphi({\mathbf w})\,,\qquad A^p{\mathbf w}=\varphi({\mathbf v})\,.
$$
\end{observation}

The following simple observation will be used as a tool in the proof of Proposition~\ref{w v}.
\begin{observation}\label{technical}
\begin{enumerate}
\item If $BA^kB$ is a factor of ${\mathbf u}_\beta$, then $k=p$ or $k=q+p$.
\item If $uB$ is a factor of ${\mathbf u}_\beta$ such that it has the prefix $A^pB$ or $A^{q+p}$, then there exists a unique factor $u'$ of ${\mathbf u}_\beta$ satisfying $\varphi(u'A)=uB$. Moreover, $|u'A|<|uB|$.
\end{enumerate}
\end{observation}

\begin{proposition}\label{w v}
Let $u$ be a~factor of ${\mathbf u}_\beta$, let $\hat{w}$ be the prefix of ${\mathbf w}$ of length $|u|$ and let $\hat{v}$ be the prefix of ${\mathbf v}$ of length $|u|$. Then $$|\hat{w}|_B\geq |u|_B\geq |\hat v|_B\,.$$
\end{proposition}
\pf
We will prove the statement by contradiction. Find the shortest factor $u$, for which the statement is not satisfied.
Then either $|\hat{w}|_B < |u|_B$ or $|\hat{v}|_B >|u|_B$.
\begin{enumerate}
\item Assume $|\hat{w}|_B < |u|_B$. By the minimality of $|u|$,  $|\hat{w}|_B+1 = |u|_B$.
According to the definition of $\mathbf w$ and since we have chosen $u$ of the minimal length, it holds, with regard to the first statement of Observation~\ref{technical},
$$
\hat w=B\cdots BA^\ell, \quad \text{where $\ell>p$, \quad and} \quad u=B\cdots B.
$$
It also follows from the first statement of Observation~\ref{technical} that 
$$
\tilde{w}=A^p\hat wA^{q+p-\ell}B=A^pB\cdots BA^qA^pB \quad\, \text{and} \,\quad \tilde{u}=A^pu
$$
are factors of ${\mathbf u}_\beta$, where the block $\hat wA^{q+p-\ell}B$ is necessarily a prefix of ${\mathbf w}$. Then the second statement of Observation~\ref{technical} implies that there exist factors $vA$ and $u'$ of ${\mathbf u}_\beta$, the factor $u'$ being shorter than $u$, such that
$$
\tilde{w}=\varphi(vA) \quad \text{and} \quad \tilde{u}=\varphi(u')\,.
$$
Furthermore, $vA$ is a~prefix of $\mathbf v$ by Observation~\ref{w v obs}.

Since $|\tilde{w}|_B = |\tilde{u}|_B$, it follows that $|vA|_A = |u'|_A$, hence $|v|_A<|u'|_A$. Since moreover $0<|\tilde{w}|-|\tilde{u}|\leq q$, one can deduce that $|vA|$ is greater than $|u'|$ by~$1$, hence $|v|=|u'|$. That means that we have found a~prefix $v$ of $\mathbf v$ and a~factor $u'$ of ${\mathbf u}_\beta$ such that $|v|=|u'|<|u|$ and $|v|_A<|u'|_A$, i.e., $|v|_B > |u'|_B$. This is a~contradiction with the assumption that $u$ is shortest possible.

\item Suppose $|\hat{v}|_B > |u|_B$. By the minimality of $|u|$,  $|\hat{v}|_B = |u|_B+1$. Using the definition of $\mathbf v$, the minimality of $|u|$, and the first statement of Observation~\ref{technical}, we deduce that
$$
\hat v=A^qA^p\cdots B \quad \text{and} \quad u=A^k B\cdots A^\ell\,,
$$
where both $k$ and $\ell$ are larger than $p$. We apply Observation~\ref{technical} once more to find that
$$
\tilde{u}=A^{q+p-k}uA^{j}B=A^qA^pB\cdots A^qA^pB \quad \text{for certain $j$, $0\leq j\leq q+p-\ell$}\,,
$$
is a~factor of ${\mathbf u}_\beta$, and that there exists a~factor $u'A$ of ${\mathbf u}_\beta$ such that $\tilde{u}=\varphi(u'A)$.
The second statement of Observation~\ref{technical} together with Observation~\ref{w v obs} imply that $\hat v=\varphi(\hat w)$, where $\hat w$ is a~prefix of $\mathbf w$, and that $\hat{w}$ is shorter than $\hat{v}$.

Since $|\hat{v}|_B = |\tilde{u}|_B$, one can see that $|\hat{w}|_A=|u'A|_A$, hence $|\hat{w}|_A>|u'|_A$. Since moreover $\tilde u$ is longer than $\hat v$, the factor $u'A$ is also longer than $\hat w$. Let us cut off a suffix of $u'A$ of length $|u'A|-|\hat w|$, and denote the resulting factor as $u''$. Then $u''$ is a factor of ${\mathbf u}_\beta$ of the same length as $\hat w$, and since we have cut off at least one $A$, it holds $|\hat w|_A>|u''|_A$, i.e., $|\hat w|_B<|u''|_B$. As the factor $u''$ is shorter than $u$ (because $|u''|=|\hat w|<|\hat v|=|u|$), we have reached a~contradiction with the minimality of $|u|$.
\end{enumerate}
\pfk

In order to determine $|\hat v|_B$ and $|\hat w|_B$ for any prefix $\hat v$ of ${\mathbf v}$ and $\hat w$ of ${\mathbf w}$,
we will need several observations. The first observation follows from the definition of $v^{(N)}$ and $w^{(N)}$ and from Observation~\ref{phi^j(A)}.
\begin{observation}\label{vwB}
For all $N\in\mathbb N_0$, it holds
\begin{gather*}
\left|v^{(N)}\right|=1+(q-1)\sum_{j=0}^{N}U_{2j}\,,\qquad \left|v^{(N)}\right|_B=(q-1)\left(1, 0\right)\sum_{j=0}^{N}M_\varphi^{2j}\left(\begin{array}{c}0\\1\end{array}\right)\,,\\
\left|w^{(N)}\right|=1+(q-1)\sum_{j=0}^{N-1}U_{2j+1}\,,\qquad \left|w^{(N)}\right|_B=1+(q-1)\left(1, 0\right)\sum_{j=0}^{N-1}M_\varphi^{2j+1}\left(\begin{array}{c}0\\1\end{array}\right)\,.
\end{gather*}
\end{observation}
As a~consequence of the form of $v^{(N)}$ and $w^{(N)}$, we obtain another observation.
\begin{observation}\label{rozdelSimple}
Let $\hat{v}$ be a~prefix of $\mathbf v$ and $\hat{w}$ a~prefix of $\mathbf w$ of length $n\in\mathbb N$.
\begin{itemize}
\item Find $M \in \mathbb N_0$ such that $|v^{(M)}|\leq n<|v^{(M+1)}|$. Then
$
\hat{v}=v^{(M)}\hat{u}\,,
$
where $\hat{u}$ is a~prefix of ${\mathbf u}_\beta$.
\item Find $N \in \mathbb N_0$ such that $|w^{(N)}|\leq n<|w^{(N+1)}|$. Then
$
\hat{w}=w^{(N)}\hat{u}'\,,
$
where $\hat{u}'$ is a prefix of ${\mathbf u}_\beta$.
\end{itemize}
\end{observation}

The last observation we need for determining $|\hat v|_B$ and $|\hat w|_B$ follows from Proposition~\ref{struktura pref. u} and Observation~\ref{phi^j(A)}.
\begin{observation}\label{u_Bsimple}
Let $\langle n\rangle_{U}=(d_{k},d_{k-1},\ldots,d_1,d_0)$. Then the prefix $u$ of ${\mathbf u}_\beta$ of length $n$ satisfies
\begin{equation*}
|u|_B=\left(1, 0\right)\left(\sum_{i=0}^{k}d_i M_\varphi^i\right)\left(\begin{array}{c}0\\1\end{array}\right)\,.
\end{equation*}
\end{observation}
Observations~\ref{vwB},~\ref{rozdelSimple}, and~\ref{u_Bsimple} lead to the formulae for $|\hat v|_B$ and~$|\hat w|_B$.
\begin{lemma}\label{lemma.v_Bsimple}
Let $\hat{v}$ be the prefix of $\mathbf v$ of length $n \in \mathbb N$. Find $M \in \mathbb N_0$ such that $|v^{(M)}|\leq n<|v^{(M+1)}|$. Then
\begin{equation}\label{v_Bsimple}
|\hat v|_B=\left(1,0\right)\left((q-1)\sum_{i=0}^{M}M_\varphi^{2i}+\sum_{i=0}^{k}d_i M_\varphi^i\right)\left(\begin{array}{c}0\\1\end{array}\right),
\end{equation}
where
$
\langle n-|v^{(M)}| \rangle_{U}=\langle n-1-(q-1)\sum_{i=0}^{M}U_{2i}\rangle_{U}=(d_{k}, d_{k-1},\ldots,d_1,d_0).
$
\end{lemma}
\begin{lemma}\label{lemma.w_Bsimple}
Let $\hat{w}$ be the prefix of $\mathbf w$ of length $n \in \mathbb N$. Find $N \in \mathbb N_0$ such that $|w^{(N)}|\leq n<|w^{(N+1)}|$. Then
\begin{equation}\label{w_Bsimple}
|\hat w|_B=1+\left(1, 0\right)\left((q-1)\sum_{i=0}^{N-1}M_\varphi^{2i+1}+\sum_{i=0}^{\ell}c_i M_\varphi^i\right)\left(\begin{array}{c}0\\1\end{array}\right),
\end{equation}
where
$
\langle n-|w^{(N)}|\rangle_{U}=\langle n-1-(q-1)\sum_{i=0}^{N-1}U_{2i+1}\rangle_{U}=(c_{\ell}, c_{\ell-1},\ldots,c_1,c_0).
$
\end{lemma}

At this moment, we could provide a~formula for the Abelian complexity using Proposition~\ref{AC=rozdil}. However, we prefer to find a~more elegant expression for $\AC(n)$. For this purpose we will use the following observation, which can be easily proven by mathematical induction.
\begin{observation}\label{sendvicSimple}
For all $j\in\mathbb N_0$, it holds
$$U_{2j}\leq |v^{(j)}|<U_{2j+1}\leq |w^{(j+1)}|<U_{2j+2}.$$
\end{observation}

The main theorem of this section follows.
\begin{theorem}\label{ACSimple}
Let ${\bf u}_\beta$ be the fixed point of the morphism $\varphi$ defined in~\eqref{simple} with $q>1$.
Let $n\in\mathbb N_0$, and let $J\in\mathbb N_0$ satisfy $U_{J}\leq n < U_{J+1}$, where $(U_J)_{J=0}^{\infty}$ is the sequence defined in~\eqref{Un}.
\begin{itemize}
\item
If $J$ is even: put $N:=\frac{J}{2}$, and put $M:=\frac{J}{2}$ if $|v^{(\frac{J}{2})}|\leq n$ and $M:=\frac{J}{2}-1$ otherwise.
\item
If $J$ is odd: put $M:=\frac{J-1}{2}$, and put $N:=\frac{J+1}{2}$ if $|w^{(\frac{J+1}{2})}|\leq n$ and $N:=\frac{J-1}{2}$ otherwise.
\end{itemize}
Then the Abelian complexity $\AC(n)$ of the infinite word ${\mathbf u}_\beta$ is given by the formula
$$
\AC(n)=2+\left(1,  0\right)\left((q-1)\left[(M_\varphi+I)^{-1}(M_\varphi^{2N}-I)-(M-N+1)M_\varphi^{2N}\right]+\sum_{i=0}^{J}(c_i-d_i) M_\varphi^i\right)\left(\begin{array}{c}0\\1\end{array}\right),
$$
where $$\begin{array}{lllll}
\langle n-|w^{(N)}|\rangle_{U}&=&\langle n-1-(q-1)\sum_{i=0}^{N-1}U_{2i+1}\rangle_{U}&=&(c_{J}, c_{J-1},\ldots,c_1,c_0)\,,\\
\langle n-|v^{(M)}| \rangle_{U}&=&\langle n-1-(q-1)\sum_{i=0}^{M}U_{2i}\rangle_{U}&=&(d_{J}, d_{J-1},\ldots,d_1,d_0)\,.\\
\end{array}$$
\end{theorem}
\pf
Due to Proposition~\ref{AC=rozdil}, it holds $\AC(n)=1+|\hat w|_B-|\hat v|_B$, where $\hat v$, $\hat w$ are the prefixes of $\mathbf v$ and $\mathbf w$, respectively, of length $n$. By Lemma~\ref{lemma.v_Bsimple}, we have
\begin{equation}\label{vvv}
|\hat v|_B=\left(1,0\right)\left((q-1)\sum_{i=0}^{M}M_\varphi^{2i}+\sum_{i=0}^{k}d_i M_\varphi^i\right)\left(\begin{array}{c}0\\1\end{array}\right),
\end{equation}
where $M$ is given by $|v^{(M)}|\leq n<|v^{(M+1)}|$ and we can put $k:=J$ because $n-|v^{(M)}|<n<U_{J+1}$.
Similarly, Lemma~\ref{lemma.w_Bsimple} gives
$$
|\hat w|_B=1+\left(1, 0\right)\left((q-1)\sum_{i=0}^{N-1}M_\varphi^{2i+1}+\sum_{i=0}^{\ell}c_i M_\varphi^i\right)\left(\begin{array}{c}0\\1\end{array}\right),
$$
where $N$ is determined by $|w^{(N)}|\leq n<|w^{(N+1)}|$ and we can set $\ell:=J$ because $n-|w^{(N)}|<n<U_{J+1}$.

Now, using the assumption $U_{J}\leq n < U_{J+1}$ together with Observation~\ref{sendvicSimple}, we find:
\begin{itemize}
\item If $J$ is even, then $N=\frac{J}{2}$, and $M=\frac{J}{2}-1$ or $M=\frac{J}{2}$;
\item if $J$ is odd, then $M=\frac{J-1}{2}$, and $N=\frac{J+1}{2}$ or $N=\frac{J-1}{2}$.
\end{itemize}
Note that it always holds $M=N-1$ or $M=N$, which allows us to rewrite the first sum in~\eqref{vvv} using the equality
$$
\sum_{i=0}^{M}M_\varphi^{2i}=\sum_{i=0}^{N-1}M_\varphi^{2i}+(M-N+1)M_\varphi^{2N}\,.
$$
Now it suffices to substitute the expressions for $|\hat v|_B$ and $|\hat w|_B$ into~\eqref{formule}, and after a~simple manipulation (note that $\sum_{i=0}^{N-1}M_\varphi^{2i+1}-\sum_{i=0}^{N-1}M_\varphi^{2i}=(M_\varphi+I)^{-1}(M_\varphi^{2N}-I)$), we obtain the sought formula for $\AC(n)$.
\pfk
\begin{example}
Let $p=3$ and $q=2$. Let us calculate $\AC(n)$ for $n=7$. Using~\eqref{Un}, we obtain $(U_n)_{n=0}^{\infty}=(1,4,14,\ldots)$.
We can see that $U_1=4 \leq 7 <U_2=14$, i.e., $J=1$ is odd in Theorem~\ref{ACSimple}, and consequently we have $M=\frac{J-1}{2}=0$. Since $|w^{(\frac{J+1}{2})}|=|w^{(1)}|=5\leq 7$, we get $N=\frac{J+1}{2}=1$. Now, we have to determine the normal $U$-representations:
$$\langle 7-|w^{(1)}|\rangle_U=\langle 2\rangle_U=(0,2)=(c_1,c_0) \quad \text{and} \quad \langle 7-|v^{(0)}|\rangle_U=\langle 5 \rangle_U=(1,1)=(d_1,d_0).$$
Finally, using Theorem~\ref{ACSimple}, we can calculate $\AC(7)=2+(1,0)\left((M_\varphi-I)+(I-M_\varphi)\right)\left(\begin{matrix}0\\1\end{matrix}\right)=2$.
Let us check the obtained value of $\AC(7)$ by determining $\hat w$ and $\hat v$, the prefixes of length $7$ of $\mathbf w$ and $\mathbf v$, respectively.
We have $\hat w=BAAABAA$ and $\hat v=AAAAABA$, consequently indeed $\AC(7)=|\hat w|_B-|\hat v|_B+1=2$.
\end{example}
\section{Conclusion}
We have determined the Abelian complexity for a~class of binary infinite words -- infinite words associated with
quadratic Parry numbers. Our method consisted in the construction of infinite words $\mathbf v$ and $\mathbf w$
whose prefixes are factors of the infinite word in question and such that the prefixes of $\mathbf v$ are the richest
in $A$'s (among all factors of the studied word ${\bf u}_\beta$) and the prefixes of $\mathbf w$ the richest in $B$'s.
Thus, besides deriving the formulae for $\AC(n)$, another important result of this paper is the description of $\mathbf v$ and $\mathbf w$ in the case of infinite words associated with quadratic simple Parry numbers (in the non-simple Parry case the words have been determined already in~\cite{BPT}).

\section{Acknowledgements}
We would like to thank Edita Pelantov\'a for her careful reading and fruitful comments.
We acknowledge financial support by the Czech Science Foundation grant 201/09/0584, by
the grants MSM6840770039 and LC06002 of the Ministry of Education, Youth, and Sports of the Czech Republic.

\end{document}